\documentclass[twoside,12pt]{article}

\usepackage{amssymb,amsmath}

\newenvironment{proof}{\begin{trivlist}\item[]{\it
Proof.}}{\hfill$\square$\end{trivlist}}

\newtheorem{theorem}{Theorem}[section]

\newtheorem{proposition}[theorem]{Proposition}

\input xy
\xyoption{matrix}\xyoption{arrow}
\def\fiddle#1{{\raise0.8ex\hbox{$#1$}}}

\def\Abar{\overline{A}}
\def\Bbar{\overline{B}}
\def\Cbar{\overline{C}}
\def\Hbar{\overline{H}}
\def\phibar{\overline{\phi}}
\def\psibar{\overline{\psi}}

\def\Atil{\widetilde{A}}
\def\Btil{\widetilde{B}}
\def\Ctil{\widetilde{C}}
\def\phitil{\widetilde{\phi}}
\def\psitil{\widetilde{\psi}}

\def\mc{{\mathbb C}}
\def\mq{{\mathbb Q}}

\def\O{{\cal O}}

\def\kc{k^\circ}
\def\co{\operatorname{co}}
\def\id{\operatorname{id}}
\def\spn{{\mathfrak{sp}}_{2n}}
\def\tr{{\rm tr}}
\def\Alt{{\rm Alt}}

\def\Oq{\O_q}
\def\OqMn{\Oq(M_n(k))}
\def\OqMmt{\Oq(M_{m,t}(k))}
\def\OqMtn{\Oq(M_{t,n}(k))}
\def\OqMmn{\Oq(M_{m,n}(k))}
\def\OqMnr{\Oq(M_{n,r}(k))}
\def\OqGLt{\Oq(GL_t(k))}
\def\OqGLn{\Oq(GL_n(k))}
\def\OqSLr{\Oq(SL_r(k))}

\def\OqMnR{\Oq(M_n(R))}
\def\OqMmtR{\Oq(M_{m,t}(R))}
\def\OqMtnR{\Oq(M_{t,n}(R))}
\def\OqMmnR{\Oq(M_{m,n}(R))}
\def\OqMnrR{\Oq(M_{n,r}(R))}
\def\OqMtR{\Oq(M_t(R))}
\def\OqGLnC{\Oq(GL_n({\mathbb C}))}
\def\OqGLtR{\Oq(GL_t(R))}
\def\OqGLnR{\Oq(GL_n(R))}
\def\OqSLrR{\Oq(SL_r(R))}
\def\OqSpn{\Oq(Sp_{2n}(k))}
\def\Oqspn{\Oq(\mathfrak{sp}\, k^{2n})}
\def\Oqspnm{\Oq((\mathfrak{sp}\, k^{2n})^m)}
\def\OqAltm{\Oq(\Alt_m(k))}
\def\Uqspn{U_q(\spn(k))}

\def\goesto{\longrightarrow}

\date{}

\begin{document}

\title{Quantized coinvariants at transcendental $q$}
\author{K R Goodearl and T H Lenagan\footnote{This research was partially
supported by NSF
research grant DMS-9622876 and by NATO Collaborative
Research Grant CRG.960250. 
}
}

\maketitle 

\begin{abstract} A general method is developed for deriving Quantum First
and Second Fundamental Theorems of Coinvariant Theory from classical
analogs in Invariant Theory, in the case that the quantization parameter $q$
is transcendental over a base field.  Several examples are given
illustrating the utility of the method; these recover earlier results of
various researchers including Domokos, Fioresi, Hacon, Rigal, Strickland,
and the present authors. 
\end{abstract}

\medskip


\section*{Introduction and background}

In the classic terminology of Hermann Weyl \cite{Weyl}, a full solution to
any invariant theory problem should incorporate a \emph{First Fundamental
Theorem}, giving a set of generators (finite, where possible) for the ring
of invariants, and a \emph{Second Fundamental Theorem}, giving
generators for the ideal of relations among the generators of the ring of
invariants. Many of the classical settings of Invariant Theory have
quantized analogs, and one seeks corresponding analogs of the
classical First and Second Fundamental Theorems. However, the setting must
be dualized before potential quantized analogs can be framed, since there
are no quantum analogs of the original objects, only quantum analogs of
their coordinate rings. Hence, one rephrases the classical results in terms
of rings of \underbar{co}invariants (see below), and then seeks quantized
versions of these. A first stumbling block is that in general these
coactions are not algebra homomorphisms, and so at the outset it is not
even obvious that the coinvariants form a subalgebra. However, this can
often be established; cf \cite{GLRrims}, Proposition 1.1;
\cite{DLjpaa}, Proposition 1.3.

Typically, a classical invariant theoretical setting is quantized uniformly
with respect to a parameter $q$, that is, there is a family of rings of
coinvariants to be determined, parametrized by nonzero scalars in a base
field, such that the case $q=1$ is the classical one. (Many authors,
however, restrict attention to special values of $q$, such as those
transcendental over the rational number field, and do not address the
general case.) As in the classical setting, one usually has identified
natural candidates to generate the ring of coinvariants. Effectively, then,
one has a parametrized inclusion of algebras (candidate subalgebras inside
the algebras of coinvariants), which is an equality at $q=1$, and one seeks
equality for other values of $q$. In the best of all worlds, the equality
at $q=1$ could be ``lifted'' by some general process to equality at
arbitrary $q$. Lifting to transcendental $q$ has been done succesfully in
some cases, by ad hoc methods -- see, for example, \cite{FiHa}, \cite{Str}. 
Also, an early version of
\cite{GLRrims} for the transcendental case was obtained in this way. 
Quantum Second Fundamental Theorems can be approached in a similar
manner. We develop here a general method for lifting equalities of
inclusions from $q=1$ to transcendental $q$, which applies to many
analyses of quantized coinvariants.
\medskip

In order to apply classical results as indicated above, we must be able to
transform invariants to coinvariants. For morphic actions of algebraic
groups, the setting of interest to us, invariants and coinvariants are
related as follows.
Suppose that $\gamma: G\times V\rightarrow V$ is a morphic action of an
algebraic group $G$ on a variety $V$. This action induces an action of $G$ on
$\O(V)$, where $(x.f)(v)= f(x^{-1}.v)$ for $x\in G$, $f\in \O(V)$, and $v\in
V$. The invariants for this action are, of course, those functions in $\O(V)$
which are constant on $G$-orbits. The comorphism of $\gamma$ is an
algebra homomorphism
$\gamma^*: \O(V) \rightarrow \O(G)\otimes \O(V)$, with respect to which
$\O(V)$ becomes a left $\O(G)$-comodule. Now a function
$f\in \O(V)$ is a coinvariant in this comodule when $\gamma^*(f)=1\otimes f$.
Since $1\otimes f$ corresponds to the function $(x,v) \mapsto f(v)$ on
$G\times V$, we see that $\gamma^*(f)=1\otimes f$ if and only if $f(x.v)=
f(v)$ for all $x\in G$ and $v\in V$, that is, if and only if $f$ is an
invariant function. To summarize:
$$\O(V)^{\co\O(G)}= \O(V)^G.$$

Quantized coordinate rings have been constructed for all complex semisimple
algebraic groups $G$. These quantized coordinate rings are Hopf algebras,
which we denote
$\Oq(G)$, since we are concentrating on the single parameter versions. In
those cases where a morphic action of $G$ on a variety $V$ has been
quantized, we have a quantized coordinate ring $\Oq(V)$ which supports an
$\Oq(G)$-coaction. The coaction is often not an algebra homomorphism, but
nonetheless -- as mentioned above -- the set of $\Oq(G)$-coinvariants in
$\Oq(V)$ is typically a subalgebra. The goal of Quantum First and Second
Fundamental Theorems for this setting is to give generators and relations
for the algebra $\Oq(V)^{\co\Oq(G)}$ of $\Oq(G)$-coinvariants in $\Oq(V)$.
We discuss several standard settings in later sections of the paper, and
outline how our general method applies. We recover known Quantum First
and Second Fundamental Theorems at transcendental $q$ in these settings,
with some simplifications to the original proofs, and in some cases
extending the range of the theorems.
\medskip

Throughout, $k$ will denote a field, which may be of arbitrary
characteristic and need not be algebraically closed.

\section{Reduction modulo $q-1$}

Throughout this section, we work with a field extension $\kc \subset k$
and a scalar $q\in k\setminus \kc$ which is transcendental over $\kc$.
Thus, the $\kc$-subalgebra $R=
\kc[q,q^{-1}] \subset k$ is a Laurent polynomial ring.  Let us denote
reduction modulo $q-1$ by overbars, that is, given any
$R$-module homomorphism $\phi:A\rightarrow B$, we write $\phibar:
\Abar\rightarrow \Bbar$ for the induced map $A/(q-1)A \rightarrow
B/(q-1)B$. 

\begin{proposition}
 Let $A \stackrel{\phi}\goesto B \stackrel{\psi}\goesto C$ be a complex of
$R$-modules, such that $C$ is torsionfree. Suppose that there are $R$-module
decompositions
$$A= \bigoplus_{j\in J} A_j, \qquad\qquad B= \bigoplus_{j\in J} B_j,
\qquad\qquad C=
\bigoplus_{j\in J} C_j$$
such that $B_j$ is finitely generated, $\phi(A_j)\subseteq B_j$, and
$\psi(B_j)\subseteq C_j$ for all $j\in J$.

If the reduced complex $\Abar \stackrel{\phibar}\goesto 
\Bbar \stackrel{\psibar}\goesto 
\Cbar$ is exact, then so is 
$$\xymatrixcolsep{3pc}
\xymatrix{
k\otimes_R A \ar[r]^{\id\otimes\phi} &k\otimes_R B
\ar[r]^{\,\id\otimes\psi\,} &k\otimes_R C. }$$
\end{proposition} 

\begin{proof} The hypotheses and the conclusions all reduce to the direct sum
components of the given decompositions, so it is enough to work in
one component. Hence, there is no loss of generality in assuming that $B$ is
finitely generated.

Let $S$ denote the localization of $R$ at the maximal ideal $(q-1)R$. 
Set 
$$\Atil= S\otimes_R A, \qquad\qquad \Btil= S\otimes_R B, \qquad\qquad
\Ctil= S\otimes_R C,$$ 
and let $\Atil \stackrel{\phitil}\goesto \Btil
\stackrel{\psitil}\goesto \Ctil$ denote the induced complex of $S$-modules. 
Since $\Atil/(q-1)\Atil$ is naturally isomorphic to $A/(q-1)A= \Abar$, and
similarly for $B$ and $C$, there is a commutative diagram

$$\xymatrixcolsep{4pc}
\xymatrix{
\fiddle{\Atil} \ar[r]^{\phitil} \ar[d]_{\alpha} &\fiddle{\Btil}
\ar[r]^{\psitil} \ar[d]_{\beta} &\fiddle{\Ctil} \ar[d]^{\gamma}\\
\fiddle{\Abar} \ar[r] ^{\phibar} &\fiddle{\Bbar} \ar[r]^{\psibar}
&\fiddle{\Cbar} }$$
\medskip

\noindent where $\alpha$, $\beta$, $\gamma$ are epimorphisms with kernels
$(q-1)\Atil$ etc.

The bottom row of the diagram is exact by hypothesis, and we claim that the
top row is exact. 
Consider an element $x\in \ker \psitil$. Chasing $x$ around the diagram, we
see that $x-\phitil(y)= (q-1)z$ for some $y\in \Atil$ and
$z\in \Btil$. Note that $(q-1)z \in \ker\psitil$. Since
$\Ctil$ is torsionfree, it follows that $z\in \ker\psitil$.
Thus $\ker\psitil \subseteq \phitil(\Atil) +(q-1)\ker\psitil$, whence
$$(\ker\psitil)/\phitil(\Atil) = (q-1) \bigl[ (\ker\psitil)/\phitil(\Atil)
\bigr].$$ 
Since $\ker\psitil$ is a finitely generated $S$-module, it follows
from Nakayama's Lemma that
$(\ker\psitil)/\phitil(\Atil) =0$, establishing the claim.

Since $k$ is a flat $S$-module, the
sequence $\xymatrixcolsep{3pc}
\xymatrix{
\fiddle{k\otimes_S \Atil} \ar[r]^{\id\otimes\phi} &\fiddle{k\otimes_S \Btil}
\ar[r]^{\,\id\otimes\psi\,} &\fiddle{k\otimes_S \Ctil} }$
is exact.  This is isomorphic to the sequence 
$$\xymatrixcolsep{3pc}
\xymatrix{
k\otimes_R A \ar[r]^{\id\otimes\phi} &k\otimes_R B
\ar[r]^{\,\id\otimes\psi\,} &k\otimes_R C }$$
 and therefore the
latter is exact.  \end{proof}

Proposition 1.1 is useful in obtaining quantized versions of both First and
Second Fundamental Theorems. For ease of application, we write out general
versions of both situations with appropriate notation, as follows.

\begin{theorem} Let $H$ be an $R$-bialgebra and $B$ a left
$H$-comodule  with structure map $\lambda: B\rightarrow H\otimes_R B$, and
let
$\phi:A\rightarrow B$ be an
$R$-module homomorphism whose image is contained in $B^{\co H}$. Assume that
$B$ and $H$ are torsionfree
$R$-modules, and suppose that there exist $R$-module
decompositions $A=
\bigoplus_{j\in J} A_j$ and $B= \bigoplus_{j\in J} B_j$ such that $B_j$
is finitely generated, $\phi(A_j)\subseteq B_j$, and $\lambda(B_j)\subseteq
H\otimes_R B_j$ for all $j\in J$.

Now $\Bbar$ is a left comodule over the $\kc$-bialgebra $\Hbar$,
and $k\otimes_R B$ is a left comodule over the $k$-bialgebra $k\otimes_R H$.
If $\Bbar^{\co \Hbar}$ equals the image of the induced map $\phibar :
\Abar \rightarrow \Bbar$, then $[k\otimes_R B]^{\co k\otimes_R H}$ equals
the image of the induced map $\id\otimes\phi: k\otimes_R A
\rightarrow k\otimes_R B$. \end{theorem}

\begin{proof} Since $R$ is a principal ideal domain, $H\otimes_R B$ is a
torsionfree $R$-module.  We may identify $\overline{H\otimes_R B}$ with
$\Hbar\otimes_{\kc} \Bbar$ and $k\otimes_R (H\otimes_R B)$ with
$(k\otimes_R H) \otimes_k (k\otimes_R B)$.  Now apply Proposition 1.1 to
the complex $A \stackrel{\,\phi\,}\goesto B \stackrel{\,\lambda
-\gamma\,}\goesto H\otimes_R B$, where $\gamma:B \rightarrow H\otimes_R
B$ is the map $b\mapsto 1\otimes b$.  \end{proof}

\begin{theorem} Let $\phi:A\rightarrow B$ be a homomorphism of
$R$-algebras, and $I$ an ideal of $A$ contained in $\ker\phi$. Assume that
$B$ is a torsionfree $R$-module, and suppose that there exist $R$-module
decompositions
$$I= \bigoplus_{j\in J} I_j, \qquad\qquad A= \bigoplus_{j\in J} A_j,
\qquad\qquad B=
\bigoplus_{j\in J} B_j$$
such that $A_j$ is finitely generated, $I_j\subseteq A_j$, and $\phi(A_j)
\subseteq B_j$ for all $j\in J$. 

If $\ker \phibar$ equals the image of $\overline{I}$, then the kernel of the
$k$-algebra homomorphism $\id\otimes\phi: k\otimes_R A \rightarrow k\otimes_R
B$ equals the image of $k\otimes_R I$. \end{theorem}

\begin{proof} Apply Proposition 1.1 to the complex
$I\stackrel{\,\eta\,}\goesto A\stackrel{\,\phi\,}\goesto 
B$, where $\eta$ is the inclusion map. \end{proof}

Proposition 1.1 can easily be adapted to yield other exactness conclusions.
In particular, suppose that, in addition to the given hypotheses, the
cokernels of $\phi$ and $\psi$ are torsionfree $R$-modules. Then one can
show that for all nonzero scalars $\lambda\in \kc$, the induced complex
$$A/(q-\lambda)A \longrightarrow B/(q-\lambda)B \longrightarrow
C/(q-\lambda)C$$
is exact. Unfortunately, in the applications to quantized coinvariants, it
appears to be a difficult task to verify the above additional hypotheses.

\section{The interior action of $GL_t$ on $M_{m,t}\times
M_{t,n}$}\label{interior} 

Let $m,n,t$ be positive integers with $t\le \min\{m,n\}$. The group $GL_t(k)$
acts on the variety $V:= M_{m,t}(k)\times M_{n,t}(k)$ via the rule
$$g\cdot(A,B)= (Ag^{-1},gB),$$
and consequently on the coordinate ring $\O(V) \cong \O(M_{m,t}(k)) \otimes
\O(M_{t,n}(k))$. The invariant theory of this action is closely related to the
matrix multiplication map
$$\mu: V \longrightarrow M_{m,n}(k)$$
and its comorphism $\mu^*$. The classical Fundamental Theorems for this
situation (\cite{DeCP}, Theorems 3.1, 3.4) are
\begin{enumerate}
\item[(1)] The ring of invariants $\O(V)^{GL_t(k)}$ equals the image of
$\mu^*$.
\item[(2)] The kernel of $\mu^*$ is the ideal of $\O(M_{m,n}(k))$ generated
by all
$(t+1)\times (t+1)$ minors.
\end{enumerate} 

As discussed in the introduction, the comorphism of the action map $\gamma:
GL_t(k) \times V
\rightarrow V$ is an algebra homomorphism $\gamma^*: \O(V) \rightarrow
\O(GL_t(k)) \otimes
\O(V)$ which makes
$\O(V)$ into a left comodule over $\O(GL_t(k))$, and the coinvariants for this
coaction equal the invariants for the action above: $\O(V)^{GL_t(k)}=
\O(V)^{\co
\O(GL_t(k))}$.

We may consider quantum versions of this situation, relative
to a parameter $q\in k^\times$, using the standard quantized coordinate rings
$\OqMmt$, $\OqMtn$, and $\OqGLt$. We write $X_{ij}$ for the standard
generators in each of these algebras; see, e.g., \cite{GLduke}, (1.1),
for the standard relations.

To quantize the classical coaction above, one combines quantized versions
of the actions of $GL_t(k)$ on $M_{m,t}(k)$ and $M_{n,t}(k)$ by right and left
multiplication, respectively. The latter quantizations are algebra morphisms
\begin{align*} 
\rho_q^* : \OqMmt &\longrightarrow \OqMmt \otimes \OqGLt \\
\lambda_q^*: \OqMtn &\longrightarrow \OqGLt\otimes \OqMtn 
\end{align*}
which are both determined by $X_{ij} \mapsto \sum_{l=1}^t X_{il}\otimes
X_{lj}$ for all $i,j$. The right comodule structure on $\OqMmt$ determines a
left comodule structure in the standard way, via the map $\tau\circ
(\id\otimes S)\circ \rho_q^*$, where $\tau$ is the flip and $S$ is the
antipode in $\OqGLt$. We set $\Oq(V)= \OqMmt\otimes \OqMtn$, which is a
left comodule over $\OqGLt$ in the standard way. The
resulting structure map $\gamma_q^*: \Oq(V)\rightarrow \OqGLt \otimes\Oq(V)$ is
determined by the rule
$$\gamma_q^*(a\otimes b)= \sum_{(a),(b)} S(a_1)b_{-1} \otimes a_0\otimes b_0$$
for $a\in \OqMmt$ and $b\in \OqMtn$, where $\rho_q^*(a)= \sum_{(a)} a_0\otimes
a_1$ and $\lambda_q^*(b)= \sum_{(b)} b_{-1}\otimes b_0$. It is, unfortunately,
more involved to work with $\gamma_q^*$ than with its classical predecessor
$\gamma^*$ because $\gamma_q^*$ is not an algebra homomorphism.

The analog of $\mu^*$ is the algebra homomorphism 
$\mu_q^*: \OqMmn \rightarrow \Oq(V)$
 such that $X_{ij} \mapsto \sum_{l=1}^t X_{il}\otimes
X_{lj}$ for all $i,j$. We can now state the Fundamental Theorems for the
quantized coinvariants in the current situation:
\begin{enumerate}
\item[$(1_q)$] The set of coinvariants $\Oq(V)^{\co\OqGLt}$ equals the image
of $\mu_q^*$ (and so is a subalgebra of $\Oq(V)$).
\item[$(2_q)$] The kernel of $\mu_q^*$ is the ideal of $\OqMmn$ generated by
all $(t+1)\times (t+1)$ quantum minors.
\end{enumerate} 

These theorems were proved for arbitrary $q$ in \cite{GLRrims}, Theorem
4.5, and
\cite{GLduke}, Proposition 2.4, respectively. 
 
To illustrate the use of Theorems 1.2 and 1.3, we specialize to the case that
$q$ is transcendental over a subfield $\kc\subset k$. Set $R= \kc[q,q^{-1}]$
as in the previous section. Since the construction of quantum matrix algebras
requires only a commutative base ring and an invertible element in that ring,
we can form
$\OqMmtR$ and $\OqMtnR$. These $R$-algebras are iterated skew polynomial
extensions of $R$, and thus are free $R$-modules, as is the algebra $\Oq(V_R)
:= \OqMmtR \otimes_R \OqMtnR$. Similarly, $\OqMtR$ is a free $R$-module and
an integral domain, and so the localization $\OqGLtR$, obtained by inverting
the quantum determinant, is a torsionfree $R$-module, as well as a Hopf
$R$-algebra. Since the restrictions of $\rho_q^*$ and $\lambda_q^*$ make
$\OqMmtR$ and $\OqMtnR$ into right and left comodules over $\OqGLtR$,
respectively, the restriction $\gamma_R^*$ of $\gamma_q^*$ makes $\Oq(V_R)$
into a left
$\OqGLtR$-comodule. Finally, $\mu_q^*$ restricts to an $R$-algebra
homomorphism $\mu_R^*: \OqMmnR \rightarrow \Oq(V_R)$. It is easily checked
that the image of
$\mu_R^*$ is contained in $\Oq(V_R)^{\co \OqGLtR}$, and that the
kernel of $\mu_R^*$ contains the ideal $I$
generated by all $(t+1)\times(t+1)$ quantum minors (cf  
\cite{GLRrims}, Proposition 2.3 and
\cite{GLduke}, (2.1)).

All of the quantum matrix algebras $\Oq(M_{\bullet,\bullet}(R))$ are
positively graded $R$-algebras, with the generators
$X_{ij}$ having degree $1$, and $\Oq(V_R)$ inherits a positive grading from
its two factors. In each of these algebras, the homogeneous components are
finitely generated free
$R$-modules. It is easily checked that
$$\mu_R^* \bigl( \OqMmnR_j \bigr) \subseteq \Oq(V_R)_{2j}
\qquad\text{and}\qquad \gamma_R^* \bigl( \Oq(V_R)_j \bigr) \subseteq \OqGLtR
\otimes_R \Oq(V_R)_j$$
for all $j\ge 0$. Moreover, the ideal $I$ of $\OqMmnR$ is homogeneous with
respect to this grading. Note that when we come to apply Theorem 1.3, we
should replace the grading on $\Oq(V_R)$ by the vector space decomposition
$\bigoplus_{j=0}^\infty \bigl( \Oq(V_R)_{2j} \oplus
\Oq(V_R)_{2j+1} \bigr)$, for instance.

The classical First and Second Fundamental Theorems (1) and (2) say that
$\overline{\mu_R^*}$ maps $\overline{\OqMmnR}$ onto the coinvariants of
$\overline{\Oq(V_R)}$, and that the kernel of $\overline{\mu_R^*}$
equals the image of $\overline{I}$. Therefore Theorems 1.2 and 1.3 yield the
quantized Fundamental Theorems $(1_q)$ and
$(2_q)$ with no further work in the transcendental case.

\section{The right action of $SL_r$ on $M_{n,r}$}

Fix positive integers $r<n$. In this section, we consider the right action
of $SL_r(k)$ on $M_{n,r}(k)$ by multiplication: $A.g= Ag$ for
$A\in M_{n,r}(k)$ and  $g\in SL_r(k)$. The First Fundamental Theorem for this
case (cf \cite{Ful}, Proposition 2, p 138) says that
\begin{enumerate}
\item[(1)] The ring of invariants $\O(M_{n,r}(k))^{SL_r(k)}$ equals the
subalgebra of $\O(M_{n,r}(k))$ generated by all $r\times r$ minors. 
\end{enumerate} 
To state the Second Fundamental Theorem for this case, let $k[X]$ be the
polynomial ring in a set of variables $X_I$, where $I$ runs over all
$r$-element subsets of $\{1,\dots,n\}$. In view of the theorem above, there is
a natural homomorphism of $\phi: k[X]\rightarrow \O(M_{n,r}(k))^{SL_r(k)}$.
\begin{enumerate}\setcounter{enumi}{1}
\item[(2)] The kernel of $\phi$ is generated by the Pl\"ucker relations
(\cite{Ful}, Proposition 2, p 138).
\end{enumerate} 
To quantize this situation, we make $\OqMnr$ into a right comodule over
$\OqSLr$ via the algebra homomorphism $\rho: \OqMnr \rightarrow \OqMnr\otimes
\OqSLr$ such that $\rho(X_{ij})= \sum_{l=1}^r X_{il}\otimes X_{lr}$ for all
$i,j$. The First Fundamental Theorem for the quantized coinvariants is the
statement
\begin{enumerate}
\item[$(1_q)$] The set of coinvariants $\OqMnr^{\co\OqSLr}$ equals the
subalgebra of $\OqMnr$ generated by all $r\times r$ quantum minors.
\end{enumerate}
This follows from work of Fioresi and Hacon under the assumptions that $k$ is
algebraically closed of characteristic $0$ and $q$ is transcendental over
some subfield of $k$ (\cite{FiHa}, Theorem 3.12) -- in fact, they prove
(in our notation) that
$(1_q)$ holds with $k$ replaced by a Laurent polynomial ring $\kc[q,q^{-1}]$.
The result is obtained for any nonzero $q$ in an arbitrary field $k$ in
\cite{KLR}, Theorem 5.2. 
The Fioresi-Hacon  version of the 
Second Fundamental Theorem yields a presentation of
$\OqMnr^{\co\OqSLr}$ in terms of generators $\lambda_{i_1\cdots i_r}$ where
$i_1\cdots i_r$ runs over all unordered sequences of $r$ distinct integers
from $\{1,\dots,n\}$. Each
$\lambda_{i_1\cdots i_r}$ is mapped to
$(-q)^{\ell(\pi)}D_{\{i_1,\dots,i_r\}}$ where $\ell(\pi)$ denotes the length
of the permutation $\pi$ sending $l\mapsto i_l$ for
$l=1,\dots,r$, and $D_{\{i_1,\dots,i_r\}}$ denotes the $r\times r$ quantum
minor with row index set $\{i_1,\dots,i_r\}$.
\begin{enumerate}
\item[$(2_q)$] The algebra $\OqMnr^{\co\OqSLr}$ is isomorphic to the quotient
of the free algebra $k\langle \lambda_{i_1\cdots i_r}\rangle$ modulo the
ideal generated by the relations (c) and (y) given in \cite{FiHa}, Theorem
3.14.
\end{enumerate}
The hypotheses on $k$ and $q$ are as above.

Now let $\kc$, $k$, $q$, $R$ be as in Section~\ref{interior}, 
and define everything over
$R$. In particular, $\OqMnrR$ is a right comodule algebra over
$\OqSLrR$. As already mentioned, Fioresi and Hacon have proved the
$R$-algebra versions of $(1_q)$ and $(2_q)$, assuming that $k$ is
algebraically closed of characteristic zero. The easier parts of their
work -- proving inclusions rather than equalities in these $R$-algebra forms
-- lead directly to the $k$-algebra forms of these theorems, as follows.

Let $R\langle \lambda_\bullet\rangle$ denote the free
$R$-algebra on generators $\lambda_{i_1\cdots i_r}$, and let $\phi_q:
R\langle \lambda_\bullet \rangle \rightarrow \OqMnrR$ be the
$R$-algebra homomorphism sending each $\lambda_{i_1\cdots i_r}$ to
$(-q)^{\ell(\pi)}D_{\{i_1,\dots,i_r\}}$. The image of $\phi_q$ is the
$R$-subalgebra generated by the $D_{\{i_1,\dots,i_r\}}$, which is contained
in $\OqMnrR^{\co \OqSLrR}$ by \cite{FiHa}, Lemma 3.5. Let $I_q$
be the ideal of $R\langle \lambda_\bullet \rangle$ generated by the relations
(c) and (y) of \cite{FiHa}, Theorem
3.14; that $I_q\subseteq \ker \phi_q$ was checked by
Fioresi
(\cite{Fio}, Proposition 2.21 and Theorem 3.6). Now $\OqMnrR$ and $R\langle
\lambda_\bullet \rangle$ are positively graded $R$-algebras, with the
$X_{ij}$ having degree $1$ and the $\lambda_\bullet$ having degree $r$. The
homogeneous components with respect to these gradings are finitely generated
free $R$-modules, $I_q$ is a homogeneous ideal of $R\langle
\lambda_\bullet \rangle$, and the map $\phi_q$ is homogeneous of degree $0$.

The classical Fundamental Theorems (1) and (2) say that $\overline{\phi_q}$
maps $\overline{R\langle
\lambda_\bullet \rangle}$ onto the coinvariants of $\overline{\OqMnrR}$, and
the kernel of $\overline{\phi_q}$ equals $\overline{I_q}$. (The latter
statement requires a change of relations using the results of \cite{Ful},
Chapter 8,
as observed in \cite{FiHa}, p 435.) Therefore Theorems 1.2 and 1.3 yield
the quantized Fundamental Theorems $(1_q)$ and $(2_q)$. Note that the
hypotheses of algebraic closure and characteristic zero on $k$ are not
needed, although $q$ is still assumed to be transcendental over a subfield
of
$k$.

\section{The right action of $Sp_{2n}$ on $M_{m,2n}$}

Fix positive integers $m$ and $n$, and consider the right action of
$Sp_{2n}(k)$ on $M_{m,2n}(k)$ by multiplication. Here we take
$Sp_{2n}(k)$ to preserve the standard alternating bilinear form on
$k^{2n}$, which we denote $\langle-,-\rangle$, and we view $M_{m,2n}(k)$
as the variety of $m$ vectors of length $2n$. Thus, each row $x_i =
(X_{i1},\dots,X_{i,2n})$ of generators in $\O(M_{m,2n}(k))$ corresponds
to the $i$-th of $m$ generic
$2n$-vectors. We shall need the functions describing the values of
$\langle-,-\rangle$ on two of these generic vectors:
$$z_{ij} := \langle x_i,x_j\rangle= \sum_{l=1}^n (X_{i,2l-1}X_{j,2l}-
X_{i,2l}X_{j,2l-1})$$
for $i,j=1,\dots,m$. 

Next, consider the variety $\Alt_m(k)$ of alternating $m\times m$
matrices over $k$; its coordinate ring is the algebra
$$\O(\Alt_m(k))= \O(M_m(k))/ \langle X_{ii},\, X_{ij}+ X_{ji} \mid i,j=
1,\dots,m \rangle.$$
Let us write $Y_{ij}$ for the coset of $X_{ij}$ in $\O(\Alt_m(k))$. There
is a morphism $\nu: M_{m,2n}(k) \rightarrow \Alt_m(k)$ given by the rule
$\nu(A)= ABA^{\tr}$, where $B$ is the matrix of the symplectic form
$\langle-,-\rangle$, and the comorphism $\nu^*: \O(\Alt_m(k)) \rightarrow
\O(M_m(k))$ sends $Y_{ij} \mapsto z_{ij}$ for all $i,j$. For any even
number $h$ of distinct indices $i_1,\dots,i_h$ from $\{1,\dots,m\}$, let
$[i_1,\dots,i_h] \in \O(\Alt_m(k))$ be the Pfaffian of the submatrix of
$(Y_{ij})$ obtained by taking the rows and columns with indices
$i_1,\dots,i_h$.

The First and Second Fundamental Theorems for the present situation
(\cite{DeCP}, Theorems 6.6, 6.7) state that
\begin{enumerate}
\item[(1)] The ring of invariants $\O(M_{m,2n}(k))^{Sp_{2n}(k)}$ equals
the subalgebra of $\O(M_{m,2n}(k))$ generated by the $z_{ij}$ for $i<j$,
that is, the image of $\nu^*$.
\item[(2)] If $m\le 2n+1$, the kernel of $\nu^*$ is zero, while if $m\ge
2n+2$, the kernel of $\nu^*$ is the ideal of $\O(\Alt_m(k))$ generated by
all the Pfaffians $[i_1,\dots,i_{2n+2}]$.
\end{enumerate}

A quantized version of this situation was studied by Strickland
\cite{Str}, who quantized the $U(\spn(k))$-module action on
$\O(M_{2n,m}(k))$ rather than its dual, the $\O(Sp_{2n}(k))$ coaction. We
transpose the matrices and indices from her paper in order to match the
notation for the classical situation used above. Since the approaches to
quantized enveloping algebras for $\spn(k)$ and quantized coordinate
rings for $Sp_n(k)$ become more complicated at roots of unity, let us
restrict our discussion to a quantum parameter $q\in k^\times$ which is
not a root of unity.

Recall that the quantized coordinate rings of the different classical
groups coact on different quantized coordinate rings of affine spaces. In
particular, $\OqSpn$ coacts on an algebra that we denote $\Oqspn$, the
{\it quantized coordinate ring of symplectic $2n$-space\/}
(\cite{RTF}, Definition 14; see \cite{Mus}, \S1.1, for a simpler set
of relations). Thus, $\O(M_{m,2n}(k))$ needs to be quantized by a
suitable algebra with $m\cdot2n$ generators, each row of which generates
a copy of $\Oqspn$. Let us call this algebra $\Oqspnm$ and take it to be
the $k$-algebra with generators $X_{ij}= x_{j,i}$ (for $i=1,\dots,m$ and
$j=1,\dots,2n$) satisfying the relations given by Strickland in
\cite{Str}, Equations (2.1), (2.2), (2.3). (The $x_{j,i}$ are
Strickland's generators for the algebra she denotes $B$; we transpose the
indices for the reasons indicated above.) Similarly, let us write
$\OqAltm$ for the $k$-algebra with generators $Y_{ij}= a_{j,i}$ for $1\le
j<i\le m$ satisfying the relations given in \cite{Str}, Equations (1.1).
As in \cite{Str}, Theorem 2.5(2), there is a $k$-algebra homomorphism
$\phi: \OqAltm \rightarrow \Oqspnm$ such that
$$\phi(Y_{ts})= \sum_{l=1}^n q^{l-1-n}X_{sl}X_{t,2n-l+1}- \sum_{l=1}^n
q^{n-l+1} X_{s,2n-l+1}X_{tl}$$
for $s<t$.

Set $B_q= \Oqspnm$ for the time being. Strickland defines an action of
$\Uqspn$ on $B_q$ (\cite{Str}, p 87). It is clear from the definition of
this action that $B_q$ is a locally finite dimensional left
$\Uqspn$-module. Hence, $B_q$ becomes a right comodule over the Hopf dual
$\Uqspn^\circ$ in the standard way (\cite{Mon}, Lemma 1.6.4(2)), and the
$\Uqspn$-invariants coincide with the $\Uqspn^\circ$-coinvariants
(\cite{Mon}, Lemma 1.7.2(2)). The standard quantized coordinate ring of
$Sp_{2n}(k)$, which we denote $\OqSpn$, is a sub-Hopf-algebra of
$\Uqspn^\circ$, and one can check that the coaction $B_q\rightarrow
B_q\otimes \Uqspn^\circ$ actually maps $B_q$ into $B_q\otimes \OqSpn$.
Consequently, $B_q$ is a right comodule over $\OqSpn$, and its
$\OqSpn$-coinvariants coincide with its $\Uqspn^\circ$-coinvariants.

Statements of First and Second Fundamental Theorems for quantized
coinvariants in the symplectic situation above can be given as follows,
where the {\it $q$-Pfaffians\/} $[i_1,\dots,i_h]$ are defined as in
\cite{Str}, p 82:
\begin{enumerate}
\item[$(1_q)$] The set of coinvariants $\Oqspnm^{\co\OqSpn}$ equals the
image of $\phi$.
\item[$(2_q)$] If $m\le 2n+1$, the kernel of $\phi$ is zero, while if
$m\ge 2n+2$, the kernel of $\phi$ is the ideal of $\OqAltm$ generated by
all the $q$-Pfaffians $[i_1,\dots,i_{2n+2}]$.
\end{enumerate}
Both statements have been proved by Strickland (modulo the changes of
notation discussed above) under the assumptions that ${\rm char}(k)=0$
and $q$ is transcendental over a subfield of $k$ (\cite{Str}, Theorem 2.5).
Via Theorems 1.2 and 1.3, we obtain the transcendental cases of $(1_q)$
and $(2_q)$ from the classical results (1) and (2) in arbitrary
characteristic. (We note that part of Strickland's development also
involves reduction modulo $q-1$. See the proof of \cite{Str}, Theorem
1.5.) As far as
we are aware, it is an open question whether $(1_q)$ and $(2_q)$ hold when
$q$ is algebraic over the prime subfield of $k$.

\section{The conjugation action of $GL_n$ on $M_n$}

Fix a positive integer $n$, and assume that $k$ has characteristic zero. 
In this section, we consider the quantum analogue of the classical
conjugation action of $GL_n(k)$ on
$M_n(k)$.  We shall need the trace
functions $\tr_i$ for $i= 1, \dots, n$, where $\tr_i$ is the sum of the
$i\times i$ principal minors.  Note that $\tr_1$ is the usual trace
function, and that $\tr_n$ is the determinant function. The First and
Second Fundamental Theorems for this situation (\cite{Kra}, Satz 3.1) can
be stated as follows:
\begin{enumerate}
\item[(1)] The ring of invariants $\O(M_n(k))^{GL_n(k)}$ equals the
subalgebra of $\O(M_n(k))$ generated by $\tr_1,\dots,\tr_n$.
\item[(2)] There are no relations among the $\tr_i$, that is,
$\tr_1,\dots,\tr_n$ are algebraically independent over $k$.
\end{enumerate} 

The coinvariants of a quantum analogue of the conjugation action have
been studied in \cite{DL}. The {\em right conjugation coaction} of
$\OqGLn$ on $\OqMn$ is the right coaction $\beta:\OqMn\goesto \OqMn
\otimes \OqGLn$ given by $\beta(u):= \sum_{(u)} u_2\otimes S(u_1)u_3$,
where we are using the Sweedler notation. In particular, $\beta(X_{ij})
= \sum_{l,m} X_{lm}\otimes S(X_{il})X_{mj}$ for all $i,j$. However, as
with the interior coaction studied in Section~\ref{interior}, the map
$\beta$ is not an algebra homomorphism. 

Recall that if $I$ and $J$ are subsets of $\{1, \dots, n\}$ of the same
size, then the quantum determinant of the quantum matrix subalgebra
generated by the $X_{ij}$ with $i\in I$ and $j\in J$ is denoted by
$[I|J]$ and called a {\em quantum minor} of the relevant quantum matrix
algebra. 
For each $i = 1, \dots, n$, define the weighted sums of principal quantum 
minors by  $\tau_i:= \sum_{I}
q^{-2w(I)}[I|I]$,  where $I$ runs through all $i$ element subsets of
$\{1, \dots, n\}$ and $w(I)$ is the sum of the
entries in the index set $I$.  
These weighted sums of principal quantum minors provide quantized
coinvariants for $\beta$ by \cite{DL}, Proposition 7.2.  One can state
First and Second Fundamental
Theorems for quantized coinvariants in this situation as follows:
\begin{enumerate}
\item[$(1_q)$] The set of coinvariants $\OqMn^{\co\OqGLn}$ is equal to 
the   subalgebra of $\OqMn$ generated by  $\tau_1, \dots, \tau_n$.
\item[$(2_q)$] The subalgebra $k[\tau_1, \dots, \tau_n]$ is a
commutative polynomial algebra of degree $n$. 
\end{enumerate}

These statements were proved for $k= \mc$ and $q$ not a root of unity in
\cite{DL}, Theorem 7.3, by using the corepresentation theory of the
cosemisimple Hopf algebra $\OqGLnC$. We give a partial extension below.

Now let $\kc = \mq$.  Suppose that $q\in k^\times$ is
transcendental over $\mq$, set $R= \mq[q,q^{-1}]$ and define
everything over $R$.  Then $\OqMnR$ is a right comodule over $\OqGLnR$
by using $\beta$ as above.  A straightforward calculation shows that the
$R$-subalgebra generated by $\tau_1, \dots, \tau_n$ is contained
in $\OqMnR^{\co \OqGLnR}$; see, for example, \cite{DL}, Proposition 7.2.
Also, note that when $q=1$ the coinvariants $\tau_1, \dots, \tau_n$
coincide with the classical traces $\tr_1, \dots, \tr_n$. Let $F(R)$
denote the free $R$-algebra on generators $\gamma_1, \dots, \gamma_n$,
and consider $F(R)$ to be graded by setting $\deg(\gamma_i) = i$. Let
$\phi_q$ be the $R$-algebra homomorphism sending $\gamma_i$ to $\tau_i$.
Then $\phi_q$ is homogeneous of degree $0$ and the image of $\phi_q$ is
contained in $\OqMnR^{\co \OqGLnR}$. The homogeneous components of
$\OqMnR$ and $F(R)$ are finitely generated free $R$-modules 
with  $\phi_q(F(R)_j) \subseteq \OqMnR_j$ and $\beta(\OqMnR_j) \subseteq 
\OqMnR_j\otimes\OqGLnR$. The classical First Fundamental Theorem (1) says
that $\overline{\phi_q}$ maps $\overline{F(R)}$ onto the coinvariants of
$\overline{\OqMnR}$. Theorem 1.2 now yields the quantized First 
Fundamental Theorem $(1_q)$ for the case that $\rm{char}(k)=0$ and $q$
is transcendental over $\mq$. 

In \cite{CW}, Corollary 2.3, Cohen and Westreich show that the $\tau_i$
commute, by exploiting the coquasitriangular structure of $\OqGLn$ (see,
e.g., \cite{Hay}, Theorem 3.1 and Proposition 4.1; a more detailed proof
for the case $k=\mc$ is given in \cite{KlSc}, Theorem 10.9).  Thus, if we
set $I$ to be the ideal of
$F(R)$ generated by the commutators
$\gamma_i\gamma_j - \gamma_j\gamma_i$, we see that $I\subseteq
\ker(\phi_q)$. It is obvious that $I$ is a homogeneous ideal. The
classical theory shows that $\ker(\overline{\phi_q}) = \overline{I}$.
Thus, Theorem 1.3 yields the Second Fundamental Theorem $(2_q)$ in the
case under discussion.

\noindent K R Goodearl,\\  
Department of Mathematics,\\ 
University of California, Santa Barbara,\\ 
CA 93106, \\
USA
\\ email: goodearl@math.ucsb.edu\\
\\
\noindent T H Lenagan,\\ 
School of Mathematics,\\ James Clerk Maxwell Building,\\
Kings Buildings, \\
Mayfield Road, \\
Edinburgh EH9 3JZ, \\
Scotland 
\\ email: tom@maths.ed.ac.uk

\end{document}